\documentclass{amsart}
\usepackage{amssymb}
\usepackage{amsfonts}
\usepackage{amsmath}
\usepackage{psfig}
\usepackage{psfrag}
\usepackage{graphicx}

\newtheorem{theorem}{Theorem}[section]

\newtheorem{proposition}[theorem]{Proposition}

\newtheorem{example}[theorem]{Example}

\newtheorem{remark}[theorem]{Remark}

\numberwithin{equation}{section}

\begin{document}

\title{Spatial chaos of Wang tiles with two symbols}

\author{Jin-Yu Chen}
\address{Department of Applied Mathematics, National Chiao Tung University, Hsinchu 30010, Taiwan}
\email{katrumichen@livemail.tw}

\author{Yu-Jie Chen}
\address{Department of Applied Mathematics, National Chiao Tung University, Hsinchu 30010, Taiwan}
\email{yujiechen0514@gmail.com}

\author{Wen-Guei Hu}
\address{College of Mathematics, Sichuan University, Chengdu, 610064, China}
\email{wghu@scu.edu.cn}

\author{Song-Sun Lin$^{\star}$}
\address{Department of Applied Mathematics, National Chiao Tung University, Hsinchu 30010, Taiwan}
\email{sslin@math.nctu.edu.tw}
\thanks{$^{\star}$The author would like to thank the National Science Council, R.O.C. (Contract No. NSC 98-2115-M-009) and
the ST Yau Center for partially supporting this research.}

\begin{abstract}
This investigation completely classifies the spatial chaos problem in plane edge coloring (Wang tiles) with two symbols. For a set of Wang tiles $\mathcal{B}$, spatial chaos occurs when the spatial entropy $h(\mathcal{B})$ is positive. $\mathcal{B}$ is
called a minimal cycle generator if
$\mathcal{P}(\mathcal{B})\neq\emptyset$ and
$\mathcal{P}(\mathcal{B}')=\emptyset$ whenever
$\mathcal{B}'\subsetneqq \mathcal{B}$, where $\mathcal{P}(\mathcal{B})$ is the set of all periodic patterns on
$\mathbb{Z}^{2}$ generated by $\mathcal{B}$.
Given a set of Wang tiles $\mathcal{B}$, write $\mathcal{B}=C_{1}\cup C_{2} \cup\cdots \cup C_{k} \cup N$, where $C_{j}$, $1\leq j\leq k$, are minimal cycle generators and $\mathcal{B}$ contains no minimal cycle generator except those contained in $C_{1}\cup C_{2} \cup\cdots \cup C_{k}$. Then, the positivity of spatial entropy $h(\mathcal{B})$ is completely determined by $C_{1}\cup C_{2} \cup\cdots \cup C_{k}$.

Furthermore, there are 39 equivalent classes of marginal positive-entropy (MPE) sets of Wang tiles and 18 equivalent classes of saturated zero-entropy (SZE) sets of Wang tiles. For a set of Wang tiles $\mathcal{B}$, $h(\mathcal{B})$ is positive if and only if $\mathcal{B}$ contains an MPE set, and $h(\mathcal{B})$ is zero if and only if $\mathcal{B}$ is a subset of an SZE set.
\end{abstract}

\maketitle

\section{Introduction}
\label{1}

The coloring of unit squares on $\mathbb{Z}^{2}$ has been studied
for many years in mathematical physics; see Baxter \cite{4,5}, Lieb \cite{13,14} and Penrose \cite{15}. In $1961$, when studying the proving theorem by
pattern recognition, Wang \cite{17} started to study the square
tiling of a plane. The unit squares with colored edges are
arranged side by side so that the adjacent tiles have the same
color; the tiles cannot be rotated or reflected. Today, such tiles
are called Wang tiles or Wang dominos \cite{6,9}.

The $2\times 2$ unit square is denoted by $\mathbb{Z}_{2\times
2}$. Let $\mathcal{S}_{p}$ be a set of $p$ $(\geq 1)$ symbols. The total set of all Wang tiles is
denoted by $\Sigma_{2\times 2}(p)$. A set
of Wang tiles $\mathcal{B}$, such that
$\mathcal{B}\subset\Sigma_{2\times 2}(p)$, is
called a basic set (of Wang tiles). Let $\Sigma(\mathcal{B})$ be
the set of all global patterns on $\mathbb{Z}^{2}$ that can be
constructed from the Wang tiles in $\mathcal{B}$, and let $\mathcal{P}(\mathcal{B})$
be the set of all periodic patterns on $\mathbb{Z}^{2}$ that can be constructed from the Wang tiles in $\mathcal{B}$. Clearly, $\mathcal{P}(\mathcal{B})\subseteq
\Sigma(\mathcal{B})$. The nonemptiness problem is to determine
whether or not $\Sigma(\mathcal{B})\neq \emptyset$.

Wang \cite{17} conjectured that any set of tiles that can
tile a plane can tile the plane periodically:
\begin{equation}\label{eqn:1.1}
\begin{array}{cccc}
\text{if} &\Sigma(\mathcal{B})\neq \emptyset,  &  \text{then} &\mathcal{P}(\mathcal{B})\neq \emptyset.
\end{array}
\end{equation}
In $1966$, using a set of $20,426$ tiles, Berger \cite{6} proved that Wang's conjecture (\ref{eqn:1.1}) is false, as the set of tiles tiled a plane aperiodically but not periodically. Subsequently, many authors reduced the number of tiles \cite{7,8,9,11,15,16}. Currently, the smallest number of tiles that can only tile a plane aperiodically is $13$, with five symbols for edge coloring (Wang tiles).

However, Chen \emph{et al.} \cite{6-0,10} proved that Wang's conjecture
holds for $p=2$ and $3$: any set of Wang tiles with three symbols that can tile
a plane can tile the plane periodically.

The most important element in proving (\ref{eqn:1.1}) in \cite{10} was the concept of the minimal cycle generator, which is introduced as follows.
A basic set $\mathcal{B}\subset\Sigma_{2\times 2}(p)$ is
called a minimal cycle generator if
$\mathcal{P}(\mathcal{B})\neq\emptyset$ and
$\mathcal{P}(\mathcal{B}')=\emptyset$ whenever
$\mathcal{B}'\subsetneqq \mathcal{B}$.
Denote the set of all minimal cycle generators by $\mathcal{C}(p)$.
Indeed, in \cite{10}, $\mathcal{C}(2)$ has $38$ members. Furthermore, under the symmetry group $D_{4}$ of
$\mathbb{Z}_{2\times 2}$ and the permutation group $S_{p}$ of
symbols of horizontal or vertical edges,
$\mathcal{C}(2)$ can be classified into six classes, as presented in Table A.1.
Then, (\ref{eqn:1.1}) is proven by showing $\Sigma(\mathcal{B})=\emptyset$ if $\mathcal{B}$ contains no minimal cycle generator.

 When (\ref{eqn:1.1}) holds, it is interesting to study the dynamical quantity of $\Sigma(\mathcal{B})$ by studying $\mathcal{P}(\mathcal{B})$. One main quantity is the spatial entropy $h(\mathcal{B})$ of $\Sigma(\mathcal{B})$. The computation of $h(\mathcal{B})$ is known to be a difficult problem for general $\mathcal{B}$. In particular, it is not easy to determine whether $h(\mathcal{B})$ is positive or not. Spatial chaos occurs when $h(\mathcal{B})>0$, and pattern formation occurs when $h(\mathcal{B})=0$ \cite{6-1}.

 This paper shows that for two symbols, the spatial chaos problem can be determined completely by using (\ref{eqn:1.1}) and minimal cycle generators. The first main theorem is as follows.

\begin{theorem}
\label{theorem:1.1}
Given $\mathcal{B}\subset\Sigma_{2\times 2}(2)$, write

\begin{equation}\label{eqn:1.2}
\mathcal{B}=C_{1}\cup C_{2} \cup\cdots \cup C_{k} \cup N,
\end{equation}
where $C_{j}$, $1\leq j\leq k$, are different minimal cycle generators and $\mathcal{B}$ contains no minimal cycle generator except those contained in $C_{1}\cup C_{2} \cup\cdots \cup C_{k}$. Then

\begin{equation}\label{eqn:1.3}
\begin{array}{ccc}
h(\mathcal{B})>0 & \text{if and only if } & h(C_{1}\cup C_{2} \cup\cdots \cup C_{k})>0.
\end{array}
\end{equation}
Furthermore, for pure cycles $C_{1}\cup C_{2} \cup\cdots \cup C_{k}$ with $1\leq k \leq 5$, 1187 equivalent classes have positive entropies and 31 equivalent classes have zero entropy.
\end{theorem}

Estimates of the lower bound of entropy for corner coloring have been made by Ban et al. \cite{3}. The method therein provides an effective means of showing $h>0$. The main element is the connecting operator, which is closely related to periodic patterns. This study establishes a similar result for Wang tiles (edge coloring). The result is presented as Theorem \ref{theorem:3.1}.

In proving (\ref{eqn:1.3}), the positivity of the entropy

\begin{equation}\label{eqn:1.4}
h(C_{1} \cup\cdots \cup C_{k})
\end{equation}
of pure cycles must firstly be determined. Equation (\ref{eqn:1.4}) is studied firstly with $k=1$ and up to $k=5$. For $k\geq 6$, entropies are positive. The symmetry group $D_{4}$ on $\mathbb{Z}_{2\times 2}$ and the permutation group $S_{2}$ greatly simplify the discussion.
The number of equivalent classes for $1\leq k\leq 5$ is $1218$. Among them, $1187$ equivalent classes have positive entropies, which can be proven by Theorem \ref{theorem:3.1}.

For the remaining $31$ equivalent classes $[C_{1} \cup\cdots \cup C_{k}]$ in Table A.2, let $\mathcal{M}(C_{1} \cup\cdots \cup C_{k})$ be the set of the maximal subsets of $\Sigma_{2\times 2}(2)$ that does not contain any minimal cycle generator except those contained in $C_{1}\cup C_{2} \cup\cdots \cup C_{k}$. For each equivalent class of the elements in $\mathcal{M}(C_{1} \cup\cdots \cup C_{k})$, after its specific properties have been studied, Propositions \ref{proposition:4.1}$\sim$\ref{proposition:4.5} introduce the upper bounds of entropy. Consequently, all of the classes have zero entropy. In particular, the spatial entropy of the 31 equivalent classes is also zero.

In contrast to edge coloring, Theorem \ref{theorem:1.1} is not valid for corner coloring, i.e., the positivity of entropy is not only determined by minimal cycle generators; see Remark \ref{remark:4.9}.

It is clear that there are $2^{16}=65,536$ sets of Wang tiles. To classify the positivity of spatial entropy of these sets, the following two concepts are useful. A set of Wang tiles $\mathcal{B}$ is called marginal positive-entropy (MPE) if $h(\mathcal{B})>0$ and $h(\mathcal{B}')=0$ whenever $\mathcal{B}'\varsubsetneqq \mathcal{B}$. A set of Wang tiles $\mathcal{B}$ is called saturated zero-entropy (SZE) if $h(\mathcal{B})=0$ and $h(\mathcal{B}')>0$ whenever $\mathcal{B}\supsetneqq \mathcal{B}$. Clearly, $h(\mathcal{B})>0$ if and only if $\mathcal{B}$ contains an MPE set of Wang tiles, and $h(\mathcal{B})=0$ if and only if $\mathcal{B}$ is a subset of an SZE set of Wang tiles.

After careful studies of results in Theorem \ref{theorem:1.1}, the second main theorem is obtained as follows.

\begin{theorem}
\label{theorem:1.2}
There are 39 equivalent classes of marginal positive-entropy sets of Wang tiles that are listed in Table A.3;
there are 18 equivalent classes of saturated zero-entropy sets of Wang tiles that are the equivalent classes of $\mathcal{B}_{M}$ in (14)$\sim$(31) in Table A.2.

\end{theorem}

The rest of this paper is organized as follows. Section 2 introduces the spatial entropy and symmetries of Wang tiles, and discusses minimal cycle generators for two symbols.
Section 3 introduces connecting operators and derives the lower bound of spatial entropy, which is important in proving the positivity of spatial entropy. Section 4 introduces some estimates of the upper bound of spatial entropy for two symbols, which imply zero entropy. Section 5 shows the main results. Some details are left in the Appendix.

\section{Entropy and symmetries of Wang tiles}
\label{2}

This section introduces the spatial entropy and symmetries of Wang tiles, and discusses minimal cycle generators for two symbols.


For given positive integers $m$ and $n$, the rectangular lattice $%
\mathbb{Z}_{m\times n}$ is defined by

\begin{equation*}
\mathbb{Z}_{m\times n}=\left\{(i,j)|0\leq i\leq m-1%
\text{ and } 0\leq j\leq n-1\right\}.
\end{equation*}
Moreover, for $(i,j)\in\mathbb{Z}^{2}$, the $2\times 2$ unit square lattice with the left-bottom vertex $(i,j)$ is defined by
\begin{equation*}
\mathbb{Z}_{2\times 2}((i,j))=\left\{(i,j),(i+1,j),(i,j+1),(i+1,j+1) \right\}.
\end{equation*}
For $m,n\geq 2$, denote the set of all local patterns with colored edges on $\mathbb{Z}_{m\times n}$ over $\mathcal{S}_{p}$, $p\geq 2$, by $\Sigma_{m\times n}(p)$.

Let $\mathcal{B}\subset\Sigma_{2\times 2}(p)$ be a set of admissible local patterns. For any $m,n\geq 2$, the set of all $\mathcal{B}$-admissible local patterns on $\mathbb{Z}_{m\times n}$ is defined by

\begin{equation}\label{eqn:2.1}
\Sigma_{m\times n}(\mathcal{B})=\left\{U\in \mathcal{S}_{p}^{\mathbb{Z}_{m\times n}}:U\mid_{\mathbb{Z}_{2\times2}((n_{1},n_{2}))}\in\mathcal{\mathcal{B}} \text{ for all }\mathbb{Z}_{2\times2}((n_{1},n_{2}))\subset\mathbb{Z}_{m\times n}   \right\}.
\end{equation}
Furthermore, let $\Sigma(\mathcal{B})$ be the set
 of all global patterns on $\mathbb{Z}^{2}$ that is generated by
$\mathcal{B}$.
Denote $\Gamma_{m\times n}(\mathcal{B})$ the cardinal number of $\Sigma_{m\times n}(\mathcal{B})$. Then, the spatial entropy $h(\mathcal{B})$ of $\mathcal{B}$ is defined by

\begin{equation}\label{eqn:2.2}
h(\mathcal{B})=\underset{m,n\rightarrow\infty}{\lim}\frac{\log \Gamma_{m\times n}(\mathcal{B})}{mn}.
\end{equation}

Now, in \cite{10}, the symmetry of the unit square $\mathbb{Z}_{2\times 2}$ is
introduced. The symmetry group of the unit square
$\mathbb{Z}_{2\times 2}$ is $D_{4}$, which is the dihedral group of order
eight. The group $D_{4}$ is generated by the rotation $\rho$,
through $\frac{\pi}{2}$, and the reflection $m$ about the
$y$-axis. Denote by
$D_{4}=\{I,\rho,\rho^{2},\rho^{3},m,m\rho,m\rho^{2},m\rho^{3}\}$.

Since, in edge coloring, the permutations of symbols in the
horizontal and vertical directions are mutually independent,
denote the permutations of symbols in the horizontal and vertical
edges by $\eta_{h}\in S_{p}$ and $\eta_{v}\in S_{p}$,
respectively. Then, for any $\mathcal{B}\subset\Sigma_{2\times
2}(p)$, define the equivalent class $[\mathcal{B}]$ of $\mathcal{B}$ by

\begin{equation}\label{eqn:2.3}
[\mathcal{B}]=\left\{\mathcal{B}'\subset\Sigma_{2\times 2}(p):\mathcal{B}'=\left(((\mathcal{B})_{\tau})_{\eta_{h}}\right)_{\eta_{v}}, \tau\in D_{4} \text{ and }\eta_{h},\eta_{v}\in S_{p} \right\}.
\end{equation}

The entropy $h(\mathcal{B})$ of $\mathcal{B}$ is clearly independent of the
choice of elements in $[\mathcal{B}]$: for any
$\mathcal{B}'\in[\mathcal{B}]$,

\begin{equation}\label{eqn:2.4}
\begin{array}{ccc}
h(\mathcal{B}')=h(\mathcal{B}).
\end{array}
\end{equation}

Next, minimal cycle generators for two symbols are discussed. For convenience of expression,
the tiles in $\Sigma_{2\times 2}(2)$ are named as follows.

\begin{equation*}
\hspace{-4.5cm}
\begin{array}{c}
\psfrag{a}{}
\psfrag{b}{ }
\includegraphics[scale=0.35]{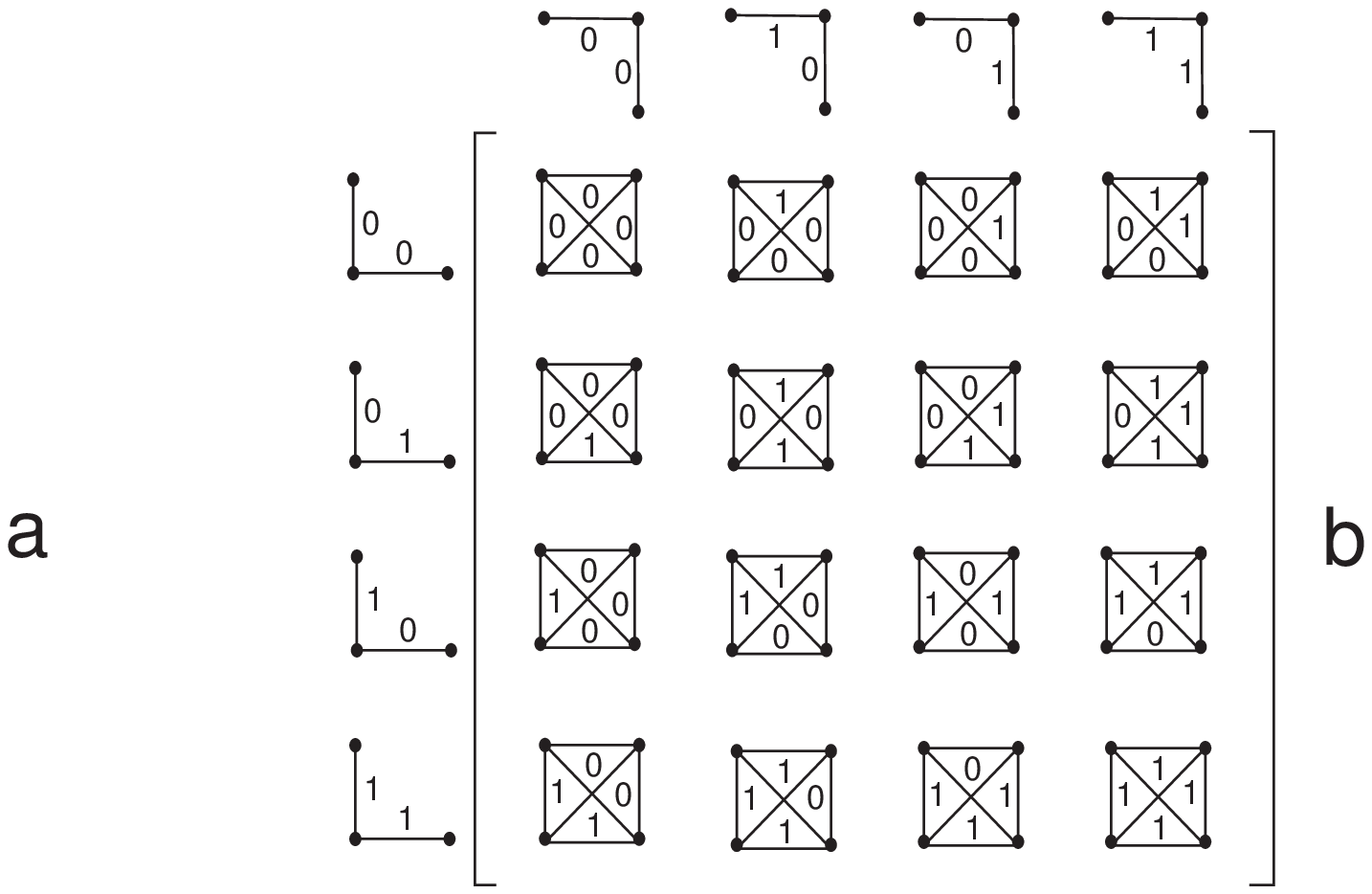}
\end{array}
\end{equation*}

\vspace{-3.0cm}
\begin{equation*}
\hspace{5.0cm}=\left[
\begin{array}{cccc}
O & E_{2} & E_{4} & R \\
E_{3} & J & B & \overline{E}_{1} \\
E_{1} & T & I & \overline{E}_{3}\\
L &  \overline{E}_{4} &  \overline{E}_{2} & E
\end{array}
\right].
\end{equation*}

\begin{equation*}
\end{equation*}
\begin{equation*}
\end{equation*}%
%
%
%
%
%
%

Now, the following theorem \cite{10} yields the six classes of $38$ minimal
cycle generators in $\mathcal{C}(2)$. Table A.1 presents the details of six equivalent classes of
$\mathcal{C}(2)$ .

\begin{theorem}
\label{theorem:2.1}
 The six equivalent classes of minimal cycle generators in
 $\mathcal{C}(2)$ are given as follows.
\begin{equation*}
\begin{array}{lll}
(1) \hspace{0.2cm} [\{O\}] & (2) \hspace{0.2cm}[\{E_{1},E_{4}\}] & (3) \hspace{0.2cm}[\{E_{1},\overline{E}_{1}\}]\\
(4) \hspace{0.2cm} [\{B,T\}]        & (5) \hspace{0.2cm}[\{E_{1},B,R\}]  & (6) \hspace{0.2cm} [\{E_{1},E_{2},B\}].
\end{array}
\end{equation*}

\end{theorem}

\section{Positive entropy}
\label{3}
This section introduces connecting operators and then derives a lower bound of spatial entropy $h(\mathcal{B})$ for $p\geq 2$.

For simplicity, only the case $p=2$ is presented for brevity. The general case can be obtained similarly.
First, the ordering matrices $\mathbf{X}_{n}$ and $\mathbf{Y}_{n}$, $n\geq 2$, are introduced to arrange systematically all local patterns in $\Sigma_{2\times n}(2)$ and $\Sigma_{n\times 2}(2)$ \cite{1,2}.

For an $n$-sequence $\overline{U}_{n}=(u_{1},u_{2},\cdots,u_{n})$ with $u_{k}\in\mathcal{S}_{2}$, $1\leq k\leq n$, $\overline{U}_{n}$ is assigned a number by applying the $n$-th order counting function $\psi\equiv\psi_{n}$:

\begin{equation}\label{eqn:3.1}
\psi(\overline{U}_{n})=\psi(u_{1},u_{2},\cdots,u_{n})=1+\underset{k=1}{\overset{n}{\sum}}u_{k}2^{(n-k)}.
\end{equation}
For $n=2$, the horizontal ordering matrix $\mathbf{X}_{2}=\underset{j=1}{\overset{4}{\sum}}\mathbf{X}_{2;j}$ and vertical ordering matrix $\mathbf{Y}_{2}=\underset{j=1}{\overset{4}{\sum}}\mathbf{Y}_{2;j}$ are defined as follows.

\begin{equation*}
\psfrag{c}{$\mathbf{X}_{2;j}=$}
\psfrag{e}{$\mathbf{Y}_{2;j}=$}
\psfrag{a}{{\footnotesize $u_{1}$}}
\psfrag{b}{{\footnotesize $u_{2}$}}
\psfrag{f}{,}
\psfrag{d}{and}
\includegraphics[scale=0.9]{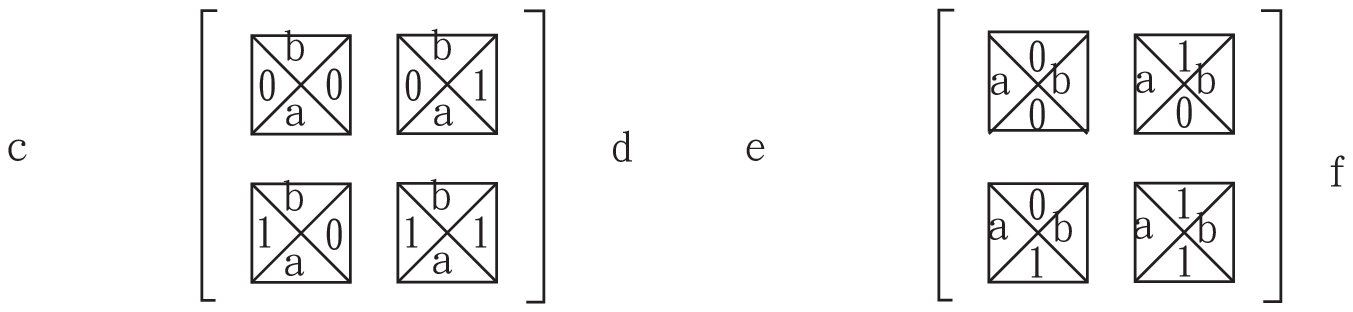}
\end{equation*}
where $u_{1},u_{2}\in\{0,1\}$ with $j=\psi(u_{1},u_{2})$.
For $n\geq 3$, the higher-order ordering matrix $\mathbf{X}_{n}=\underset{j=1}{\overset{4}{\sum}}\mathbf{X}_{n;j}$  can be recursively defined recursively by
\begin{equation}\label{eqn:3.2}
\left\{
\begin{array}{l}
\mathbf{X}_{n;1}=\left(\mathbf{X}_{2;1}\otimes \mathbf{X}_{n-1;1}\right)+\left(\mathbf{X}_{2;2}\otimes \mathbf{X}_{n-1;3}\right), \\
\mathbf{X}_{n;2}=\left(\mathbf{X}_{2;1}\otimes \mathbf{X}_{n-1;2}\right)+\left(\mathbf{X}_{2;2}\otimes \mathbf{X}_{n-1;4}\right), \\
\mathbf{X}_{n;3}=\left(\mathbf{X}_{2;3}\otimes \mathbf{X}_{n-1;1}\right)+\left(\mathbf{X}_{2;4}\otimes \mathbf{X}_{n-1;3}\right), \\
\mathbf{X}_{n;4}=\left(\mathbf{X}_{2;3}\otimes \mathbf{X}_{n-1;2}\right)+\left(\mathbf{X}_{2;4}\otimes \mathbf{X}_{n-1;4}\right),
\end{array}
\right.
\end{equation}
where $\otimes$ is the Kronecker product: if $A=[a_{i,j}]$ and $B=[b_{i,j}]$, then $A\otimes B=[a_{i,j} B]$.
Similarly, $\mathbf{Y}_{n}=\underset{j=1}{\overset{4}{\sum}}\mathbf{Y}_{n;j}$, $n\geq 3$, can be defined recursively, as in (\ref{eqn:3.2}).

For $m\geq 1$ and $1\leq \alpha\leq 4$, the horizontal connecting ordering matrix \newline $\mathbf{S}_{m;\alpha}=\left[\left(\mathbf{S}_{m;\alpha}\right)_{k,l}\right]_{2^{m}\times 2^{m}}$ is defined, where $\left(\mathbf{S}_{m;\alpha}\right)_{k,l}$ is the set of all local patterns of the form,

\begin{equation}\label{eqn:3.3}
\begin{array}{c}
\psfrag{a}{{\footnotesize $k_{1}$}}
\psfrag{b}{{\footnotesize $k_{2}$}}
\psfrag{c}{{\footnotesize $k_{m}$}}
\psfrag{k}{{\footnotesize $l_{1}$}}
\psfrag{l}{{\footnotesize $l_{2}$}}
\psfrag{m}{{\footnotesize $l_{m}$}}
\psfrag{d}{{\footnotesize $\alpha_{1}$}}
\psfrag{h}{{\footnotesize $\alpha_{2}$}}
\psfrag{o}{$\cdots$}
\psfrag{e}{{\footnotesize $b_{1}$}}
\psfrag{f}{{\footnotesize $b_{2}$}}
\psfrag{g}{{\tiny $b_{m-1}$}}
\includegraphics[scale=1.4]{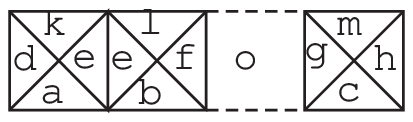}
\end{array}
\end{equation}
with $b_{j}\in\{0,1\}$, $1\leq j\leq m-1$, $\alpha=\psi(\alpha_{1},\alpha_{2})$, $k=\psi(k_{1},\ldots,k_{m})$ and $l=\psi(l_{1},\ldots,l_{m})$. Indeed,

\begin{equation}\label{eqn:3.4}
\mathbf{S}_{m;\alpha}=\mathbf{Y}_{m+1;\alpha}
\end{equation}
for all $m\geq 1$, $1\leq\alpha\leq 4$.
Notably, from (\ref{eqn:3.3}), for $m\geq 1$, each local pattern in $\mathbf{S}_{m;\alpha}$, $\alpha\in\{1,4\}$, has horizontal periodicity with period $m$, such that $\alpha_{1}=\alpha_{2}$.
Similarly, the vertical connecting ordering matrix $\mathbf{W}_{m;\alpha}=\left[\left(\mathbf{W}_{m;\alpha}\right)_{k,l}\right]_{2^{m}\times 2^{m}}$ can be defined.

Now, given a basic set $\mathcal{B}\subset \Sigma_{2\times2}(2)$, the horizontal transition matrix $\mathbb{H}_{2}=\underset{j=1}{\overset{4}{\sum}}H_{2;j}$ is defined as follows. Let $H_{2;j}=\left[h_{2;j;s,t}\right]_{2\times 2}$, where

\begin{equation}\label{eqn:3.5}
\left\{
\begin{array}{rl}
h_{2;j;s,t}=1 & \text{if }\left(\mathbf{X}_{2;j} \right)_{s,t}\in\mathcal{B}, \\
=0 & \text{otherwise. }
\end{array}
\right.
\end{equation}
As in (\ref{eqn:3.2}), for $n\geq 3$, $\mathbb{H}_{n}=\underset{j=1}{\overset{4}{\sum}}H_{n;j}$ can be defined as

\begin{equation}\label{eqn:3.6}
\left\{
\begin{array}{l}
H_{n;1}=\left(H_{2;1}\otimes H_{n-1;1}\right)+\left(H_{2;2}\otimes H_{n-1;3}\right), \\
H_{n;2}=\left(H_{2;1}\otimes H_{n-1;2}\right)+\left(H_{2;2}\otimes H_{n-1;4}\right), \\
H_{n;3}=\left(H_{2;3}\otimes H_{n-1;1}\right)+\left(H_{2;4}\otimes H_{n-1;3}\right), \\
H_{n;4}=\left(H_{2;3}\otimes H_{n-1;2}\right)+\left(H_{2;4}\otimes H_{n-1;4}\right).
\end{array}
\right.
\end{equation}
Similarly, $\mathbb{V}_{n}=\underset{j=1}{\overset{4}{\sum}}V_{n;j}$, $n\geq 3$, can be defined recursively.

Given $\mathcal{B}\subset \Sigma_{2\times2}(2)$, for $m\geq 1$ and $1\leq \alpha\leq 4$, the horizontal connecting operator $S_{m;\alpha}=\left[\left(S_{m;\alpha}\right)_{k,l}\right]_{2^{m}\times 2^{m}}$ can be defined, where
$\left(S_{m;\alpha}\right)_{k,l}$ is the cardinal number of all $\mathcal{B}$-admissible local patterns in
$\left(\mathbf{S}_{m;\alpha}\right)_{k,l}$.
Furthermore, from (\ref{eqn:3.4}),

\begin{equation}\label{eqn:3.7}
S_{m;\alpha}=V_{m+1;\alpha}
\end{equation}
for $m\geq 1$ and $1\leq \alpha\leq 4$.
Similarly, the vertical connecting operator is denoted by $\mathbb{W}_{m;\alpha}=\left[\left(W_{m;\alpha}\right)_{k,l}\right]_{2^{m}\times 2^{m}}$.

Connecting operators $S_{m;\alpha}$ or $W_{m;\alpha}$, $\alpha\in\{1,4\}$, are used to obtain a lower bound of spatial entropy $h(\mathcal{B})$, which is very useful in proving the positivity of entropy, from the following theorem.

\begin{theorem}
\label{theorem:3.1}
Given $\mathcal{B}\subset\Sigma_{2\times 2}(2)$, let $\beta_{1},\beta_{2},\cdots,\beta_{k}\in\{1,4\}$, $k\geq 1$. Then, for any $m\geq 1$,

\begin{equation}\label{eqn:3.8}
h(\mathcal{B})\geq \frac{1}{mk}\log\rho\left(S_{m;\beta_{1}}S_{m;\beta_{2}}\cdots S_{m;\beta_{k}}\right)
\end{equation}
and
\begin{equation}\label{eqn:3.9}
h(\mathcal{B})\geq \frac{1}{mk}\log\rho\left(W_{m;\beta_{1}}W_{m;\beta_{2}}\cdots W_{m;\beta_{k}}\right)
\end{equation}
where $\rho(A)$ is the maximum eigenvalue of the matrix $A$.
\end{theorem}
\textit{Proof}
For simplicity, only (\ref{eqn:3.8}) is proven. The proof of (\ref{eqn:3.9}) is similar.

For fixed $m,k\geq 1$, from the construction of $S_{m;\alpha}$, it is easy to see that for any $s\geq 1$,

\begin{equation*}
\Gamma_{(m+1)\times(sk+1)}(\mathcal{B})=\left| \mathbb{V}_{m+1}^{sk}(\mathcal{B})\right|\geq \left| \left(S_{m;\beta_{1}}S_{m;\beta_{2}}\cdots S_{m;\beta_{k}}    \right)^{s}\right|,
\end{equation*}
where $|A|$ is the sum of all entries in matrix $A$.

Since $\beta_{1},\beta_{2},\cdots,\beta_{k}\in\{1,4\}$, the local patterns in $\mathbf{S}_{m;\beta_{1}},\mathbf{S}_{m;\beta_{2}},\cdots,\mathbf{S}_{m;\beta_{k}}$ have horizontal periodicity with period $m$. Hence, for any $s,t\geq 1$,

\begin{equation*}
\Gamma_{(tm+1)\times(sk+1)}(\mathcal{B})\geq \left| \left(S_{m;\beta_{1}}S_{m;\beta_{2}}\cdots S_{m;\beta_{k}}    \right)^{s}\right|^{t}.
\end{equation*}
Therefore, by the Perron-Frobenius theorem,
\begin{equation*}
\begin{array}{rl}
h(\mathcal{B})= & \underset{s,t\rightarrow\infty}{\limsup}\frac{1}{(tm+1)(sk+1)}\log \Gamma_{(tm+1)\times(sk+1)}(\mathcal{B}) \\
\geq& \underset{s,t\rightarrow\infty}{\limsup}\frac{1}{(tm+1)(sk+1)}\log \left| \left(S_{m;\beta_{1}}S_{m;\beta_{2}}\cdots S_{m;\beta_{k}}    \right)^{s}\right|^{t}\\
=& \frac{1}{mk}\log\rho\left(S_{m;\beta_{1}}S_{m;\beta_{2}}\cdots S_{m;\beta_{k}}\right).
\end{array}
\end{equation*}
The proof is complete.
$\Box$
\medbreak
Notably, Theorem \ref{theorem:3.1} also holds for $p\geq 3$.

The following six-vertex model (or ice-type model) is used to illustrate Theorem \ref{theorem:3.1}.

\begin{example}
\label{example:3.2}
The rule of the six-vertex model is that the number of arrows that point inwards at each vertex is two, such that
 \begin{equation*}
 \mathcal{B}'_{6}=\left\{
\begin{array}{cccccc}
\includegraphics[scale=0.45]{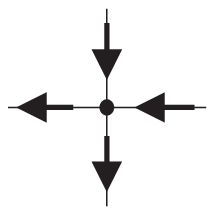}, &
\includegraphics[scale=0.45]{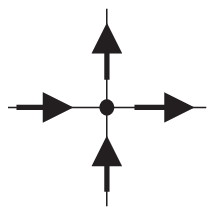}, &
\includegraphics[scale=0.45]{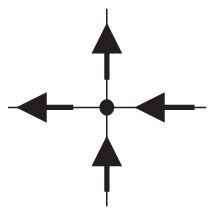}, &
\includegraphics[scale=0.45]{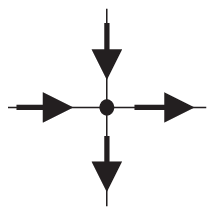}, &
\includegraphics[scale=0.45]{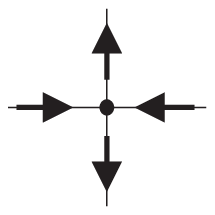}, &
\includegraphics[scale=0.45]{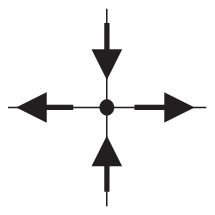}
\end{array}
\right\}.
\end{equation*}
The six-vertex model can be transformed into an edge-coloring problem with the symbols in $\mathcal{S}_{2}=\{0,1\}$, as follows:
the rightward and upward arrows in each pattern in $\mathcal{B}'_{6}$ are replaced by the digit (color) $1$ and the leftward and downward arrows in each pattern are replaced by the digit (color) $0$.  Indeed,
\begin{equation*}
 \mathcal{B}_{6}=\left\{
\begin{array}{cccccc}
\includegraphics[scale=0.5]{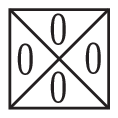}, &
\includegraphics[scale=0.5]{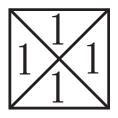}, &
\includegraphics[scale=0.5]{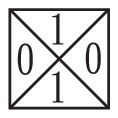}, &
\includegraphics[scale=0.5]{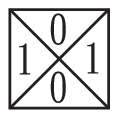}, &
\includegraphics[scale=0.5]{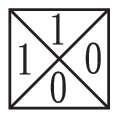}, &
\includegraphics[scale=0.5]{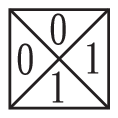}
\end{array}
\right\}.
\end{equation*}
In \cite{14}, that $h(\mathcal{B}_{6})=\frac{3}{2}\log \left(\frac{4}{3}\right)$ was proven.

Clearly,
\begin{equation*}
\begin{array}{cccc}
V_{2;1}=\left[
\begin{array}{cc}
1 & 0 \\
0 & 1
\end{array}
\right],
&
V_{2;2}=\left[
\begin{array}{cc}
0 & 0 \\
1 & 0
\end{array}
\right],
&
V_{2;3}=\left[
\begin{array}{cc}
0 & 1 \\
0 & 0
\end{array}
\right],
&
V_{2;4}=\left[
\begin{array}{cc}
1 & 0 \\
0 & 1
\end{array}
\right].
\end{array}
\end{equation*}
Then,
\begin{equation*}
S_{2;1}=V_{3;1}=\left(V_{2;1}\otimes V_{2;1}\right)+\left(V_{2;2}\otimes V_{2;3}\right)=
\left[
\begin{array}{cccc}
1 & 0 & 0 & 0 \\
0 & 1 & 0 & 0 \\
0 & 1 & 1 & 0 \\
0 & 0 & 0 & 1
\end{array}
\right]
\end{equation*}
and
\begin{equation*}
S_{2;4}=V_{3;4}=\left(V_{2;3}\otimes V_{2;2}\right)+\left(V_{2;4}\otimes V_{2;4}\right)=
\left[
\begin{array}{cccc}
1 & 0 & 0 & 0 \\
0 & 1 & 1 & 0 \\
1 & 0 & 1 & 0 \\
0 & 0 & 0 & 1
\end{array}
\right].
\end{equation*}
That $\rho(S_{2;1}S_{2;4})=\frac{3+\sqrt{5}}{2}$ can be easily shown. Therefore, from (\ref{eqn:3.8}), $h(\mathcal{B}_{6})\geq \frac{1}{4}\log \left( \frac{3+\sqrt{5}}{2}\right)$.

\end{example}

Notably, the well-known eight-vertex model also can be transformed into
 \begin{equation*}
 \mathcal{B}_{8}=\left\{
\begin{array}{cccccccc}
\includegraphics[scale=0.5]{t61.eps}, &
\includegraphics[scale=0.5]{t62.eps}, &
\includegraphics[scale=0.5]{t63.eps}, &
\includegraphics[scale=0.5]{t64.eps}, &
\includegraphics[scale=0.5]{t65.eps}, &
\includegraphics[scale=0.5]{t66.eps}, &
\includegraphics[scale=0.5]{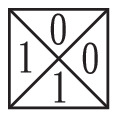}, &
\includegraphics[scale=0.5]{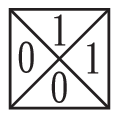}
\end{array}
\right\},
\end{equation*}
which contains $\mathcal{B}_{6}$ and then it also has positive entropy. Furthermore, it can be easily proven that $h(\mathcal{B}_{8})=\log 2$.

In fact, in the application of (\ref{eqn:3.8}) or (\ref{eqn:3.9}) to prove that the entropy of the 1187 equivalent classes that are not presented in Table A.2 is positive, the maximal value of $m$ in Theorem \ref{theorem:3.1} is five. However, in applying Theorem \ref{theorem:3.1} to test the remaining 31 equivalent classes that are presented in Table A.2, the maximum eigenvalues $\rho\left(S_{m;\beta_{1}}S_{m;\beta_{2}}\cdots S_{m;\beta_{k}}\right)$ and $\rho\left(W_{m;\beta_{1}}W_{m;\beta_{2}}\cdots W_{m;\beta_{k}}\right)$ still equal to one even when $m\geq 10$. Therefore, it is reasonable to guess that the entropies of them are zero, which will be shown in next section.

An example of the 1187 equivalent classes with $m=5$ follows.

\begin{example}
\label{example:3.3}
Consider $\mathcal{B}=\left\{
\begin{array}{cccc}
\includegraphics[scale=0.5]{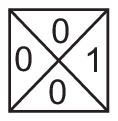},& \includegraphics[scale=0.5]{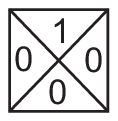},&
\includegraphics[scale=0.5]{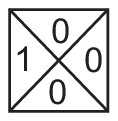}, &\includegraphics[scale=0.5]{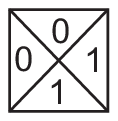}
\end{array}
\right\}$ which is the eighth case of Table A.3. Clearly,

\begin{equation*}
\begin{array}{cccc}
V_{2;1}=\left[
\begin{array}{cc}
0 & 1 \\
0 & 0
\end{array}
\right],
&
V_{2;2}=\left[
\begin{array}{cc}
1 & 0 \\
1 & 0
\end{array}
\right],
&
V_{2;3}=\left[
\begin{array}{cc}
1 & 0 \\
0 & 0
\end{array}
\right],
&
V_{2;4}=\left[
\begin{array}{cc}
0 & 0 \\
0 & 0
\end{array}
\right].
\end{array}
\end{equation*}
It can be verified that
\begin{equation*}
\rho(S_{5;1}S_{5;1}S_{5;4})=2+\sqrt{3}.
\end{equation*}

Therefore, from (\ref{eqn:3.8}),  $h(\mathcal{B})\geq\frac{1}{15}\log(2+\sqrt{3}) >0 $.

\end{example}

From (\ref{eqn:3.8}) and (\ref{eqn:3.9}), for $\mathcal{B}=C_{1}\cup C_{2} \cup\cdots \cup C_{k}\subset\Sigma_{2\times 2}(2)$, $1\leq k\leq 5$, carefully checking the spatial entropies $h(\mathcal{B})$ of the 1218 equivalent classes yields the following result.

Before showing the following theorem, for simplifying the discussion of the positivity problem, the following notation is introduced.
For two equivalent classes $[\mathcal{B}_{1}]$ and $[\mathcal{B}_{2}]$, if there exist $\mathcal{B}'\in[\mathcal{B}_{1}]$ and $\mathcal{B}''\in[\mathcal{B}_{2}]$ such that $\mathcal{B}'\supseteq\mathcal{B}''$, then denote $[\mathcal{B}_{1}]\supseteq [\mathcal{B}_{2}]$. Clearly, if $[\mathcal{B}_{1}]\supseteq [\mathcal{B}_{2}]$ with $h([\mathcal{B}_{2}])>0$, then $h([\mathcal{B}_{1}])>0$; if $[\mathcal{B}_{1}]\supseteq [\mathcal{B}_{2}]$ with $h([\mathcal{B}_{1}])=0$, then $h([\mathcal{B}_{2}])=0$.

\begin{theorem}
\label{theorem:3.4}
Given $\mathcal{B}\subset\Sigma_{2\times 2}(2)$, write
$\mathcal{B}=C_{1}\cup C_{2} \cup\cdots \cup C_{k}$
for $k\geq 1$. If $1\leq k \leq 5$ and $\mathcal{B}$ is not in the 31 equivalent classes that are presented in Table A.2, then $h(\mathcal{B})>0$. Furthermore, if $k=5$, then $h(\mathcal{B})>0$.
\end{theorem}
\textit{Proof}
Consider $\mathcal{B}=C_{1}\cup C_{2} \cup\cdots \cup C_{k}$ for $1\leq k \leq 5$. From (\ref{eqn:2.3}), it can be verified that the number of equivalent classes is 1218. In the following, it will be shown that the 1187 equivalent classes $[\mathcal{B}]$ that are not in Table A.2 have positive entropies.

The 1187 equivalent classes can be simplified as follows. If $[\mathcal{B}_{1}]$ and $[\mathcal{B}_{2}]$ are two of the 1187 equivalent classes with $[\mathcal{B}_{1}]\supseteq [\mathcal{B}_{2}]$, then delete $[\mathcal{B}_{1}]$ from the 1187 equivalent classes. By a computer program, it can be checked that the 1189 equivalent classes can be reduced to the 39 equivalent classes that are presented in Table A.3.
Therefore, it suffices to prove that the 39 equivalent classes in Table A.3 have positive entropies.

For each $\mathcal{B}$ in the 39 equivalent classes $[\mathcal{B}]$ in Table A.3, (\ref{eqn:3.8}) (or (\ref{eqn:3.9})) with some $1\leq m,k\leq 5$ yields $h(\mathcal{B})>0$. The relevant details are also presented in Table A.3.
The proof is complete.
$\Box$
\medbreak

\section{Zero entropy}
\label{4}
In this section, a series of propositions are designed for estimating the upper bound of spatial entropy for the elements in $\mathcal{M}(C_{1} \cup\cdots \cup C_{k})$ of the 31 equivalent classes in Table A.2.

First, relevant notation must be specified. For $B\in\mathcal{B}\subset \Sigma_{2\times 2}(2)$, let
\newline
$b(B), t(B), l(B), r(B)\in\{0,1\}$ be the symbols of the bottom, top, left and right edges in $B$, respectively.

\begin{proposition}
\label{proposition:4.1}
Given $\mathcal{B}\subset\Sigma_{2\times 2}(2)$. If $\left(b(B_{1}),l(B_{1})\right)\neq \left(b(B_{2}),l(B_{2})\right)$ for any $B_{1},B_{2}\in\mathcal{B}$ with $B_{1}\neq B_{2}$, then $h(\mathcal{B})=0$.
\end{proposition}
\textit{Proof}
 Based on the assumption, a pattern $U_{m\times n}\in\Sigma_{m\times n}(\mathcal{B})$ is easily seen to be uniquely determined by $b\left(U_{m\times n}\mid_{\mathbb{Z}_{2\times 2}((i,0))}\right)$ and $l\left(U_{m\times n}\mid_{\mathbb{Z}_{2\times 2}((0,j))}\right)$, $0\leq i\leq m-2$ and $0\leq j\leq n-2$.

Then, $\Gamma_{m\times n}(\mathcal{B})\leq 2^{(m-1)+(n-1)}$. Therefore, $h(\mathcal{B})=0$ follows.
$\Box$
\medbreak

\begin{proposition}
\label{proposition:4.2}
Given $\mathcal{B}\subset\Sigma_{2\times 2}(2)$. If every $B\in\mathcal{B}$ with $B\neq R = $\hspace{0.1cm}\includegraphics[scale=0.4]{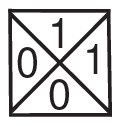} satisfies that
$\left(b(B),t(B)\right),\left(r(B),l(B)\right)\in\{(0,0),(1,0),(1,1)\}$, then $h(\mathcal{B})=0$.
\end{proposition}
\textit{Proof}
Since $\left(l(B),r(B)\right)\in\{(0,0),(0,1),(1,1)\}$ for all $B\in\mathcal{B}$, at most $n+2$ patterns on
the lattice
\begin{equation*}
\begin{array}{c}
\psfrag{a}{$\cdots$}
\psfrag{n}{ $n$\text{ copies}}
\includegraphics[scale=0.8]{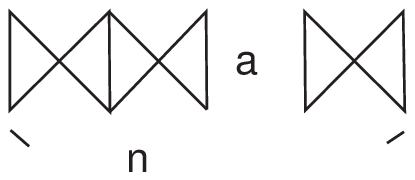}
\end{array}
\end{equation*}
can be generated by $\left(l(B),r(B)\right)$ of $B\in\mathcal{B}$.

Notably, tile $R$ is the only tile in $\mathcal{B}$ such that the symbol on bottom edge is $0$ and the symbol on top edge is $1$; this property is important in vertical tiling. However, every pattern in $\Sigma_{(n+1)\times 2}(\mathcal{B})$ can contain no more than one tile $R$. Then, every pattern in $\Sigma_{(n+1)\times (n+1)}(\mathcal{B})$ contains no more than $n$ tiles $R$.
That $\Gamma_{(n+1)\times (n+1)}(\mathcal{B})\leq  (n+2)^{3n}$ can be verified. Therefore, $h(\mathcal{B})=0$.
$\Box$
\medbreak

\begin{proposition}
\label{proposition:4.3}
Given $\mathcal{B}\subset\Sigma_{2\times 2}(2)$. If every $B\in\mathcal{B}$ satisfies that
$\left(b(B),l(B)\right)\neq (0,1)$ and $\left(t(B),r(B)\right)\neq (1,0)$, then $h(\mathcal{B})=0$.
\end{proposition}
\textit{Proof}
Clearly, the numbers (which the symbols $0$ and $1$ are regarded as) of the pattern $U_{(n+1)\times (n+1)}\in\Sigma_{(n+1)\times (n+1)}(\mathcal{B})$ along $\frac{3}{4}\pi$ lines must be decreasing from the bottom-right to the upper-left. For example,
\begin{equation*}
\begin{array}{c}
\psfrag{a}{.}
\includegraphics[scale=0.75]{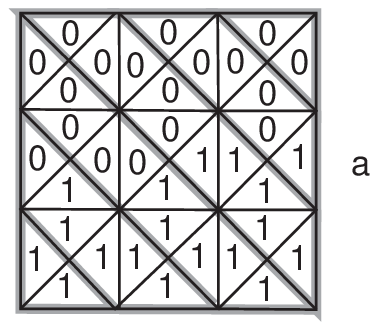}
\end{array}
\end{equation*}
Then, that $\Gamma_{(n+1) \times (n+1)}(\mathcal{B})\leq   (2n+1)^{2n}$ can be shown. Therefore, $h(\mathcal{B})= 0$.
$\Box$
\medbreak

\begin{proposition}
\label{proposition:4.4}
Given $\mathcal{B}\subset\Sigma_{2\times 2}(2)$. If $\mathcal{B}$ satisfies that

\begin{enumerate}
\item[(i)] for any $B\in\mathcal{B}$, $(l(B),r(B))\in\{(0,0),(0,1),(1,1)\}$,
\item[(ii)] if $B_{k}\in\mathcal{B}$ with $\left(l(B_{k}),r(B_{k})\right)=(0,0)$, $k\in\{1,2\}$, then
$b(B_{1})\neq b(B_{2})$,
\item[(iii)] if $B_{k}\in\mathcal{B}$ with $\left(l(B_{k}),r(B_{k})\right)=(1,1)$, $k\in\{1,2\}$, then
$b(B_{1})\neq b(B_{2})$,
\end{enumerate}

then $h(\mathcal{B})=0$.
\end{proposition}
\textit{Proof}
For $n\geq 1$, an upper bound of $\Gamma_{(n+1)\times (n+1)}(\mathcal{B})$ is obtained by the following four steps.

\begin{enumerate}
\item[Step (1):] color all horizontal edges on $\mathbb{Z}_{(n+1)\times (n+1)}$ according to condition (i). Then, the number of possible patterns is less than $(n+2)^{n}$. Notably, for any row of these patterns, no more than one such pattern \includegraphics[scale=0.4]{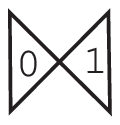} exists.

\item[Step (2):] color the vertical edges of \includegraphics[scale=0.4]{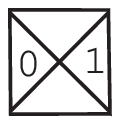} of the patterns that were obtained in step (1) by $0$ or $1$. Then, the number of possible patterns is less than $(n+2)^{n}\cdot 2^{2n}$.
\item[Step (3):] color the remaining vertical edges on the bottom row of the patterns that were obtained in step (2) by $0$ or $1$. Then, the number of the possible patterns is less than $(n+2)^{n}\cdot 2^{2n}\cdot 2^{2n}$.

\item[Step (4):] color all remaining vertical edges of the patterns that were obtained in step (3). Conditions (ii) and (iii) uniquely determine the coloring in this step.
\end{enumerate}
Then, $\Gamma_{(n+1)\times (n+1)}(\mathcal{B})\leq 2^{4n}(n+2)^{n} $ for $n\geq 1$. Therefore, $h(\mathcal{B})=0$.
$\Box$
\medbreak

%

\begin{proposition}
\label{proposition:4.5}
If $\mathcal{B}\in\left[ O,E_{1},\overline{E}_{2},R \right]$, then $h(\mathcal{B})=0$.
\end{proposition}
\textit{Proof}
First, the following three facts can be easily established.
\begin{enumerate}
\item[Fact (1):]
For tile $\overline{E}_{2}$, if a pattern $U_{4\times 4}\in \Sigma_{4\times 4}(\mathcal{B})$ with $U_{4\times 4}\mid_{\mathbb{Z}_{2\times 2}((1,1))}=\overline{E}_{2}$ exists, then $U_{4\times 4}\mid_{\mathbb{Z}_{2\times 2}((0,2))}=\overline{E}_{2}$.

\item[Fact (2):]
For tile $R$, if a pattern $U_{3\times 3}\in \Sigma_{3\times 3}(\mathcal{B})$ with $U_{3\times 3}\mid_{\mathbb{Z}_{2\times 2}((1,0))}=R$ exists, then \newline $U_{3\times 3}\mid_{\mathbb{Z}_{2\times 2}((0,1))}=R$.

\item[Fact (3):]
For tile $E_{1}$, if a pattern $U_{4\times 4}\in \Sigma_{4\times 4}(\mathcal{B})$ with $U_{4\times 4}\mid_{\mathbb{Z}_{2\times 2}((2,1))}=E_{1}$ exists, then $U_{4\times 4}\mid_{\mathbb{Z}_{2\times 2}((1,2))}$ is $E_{1}$ or $\overline{E}_{2}$.
 \end{enumerate}

 Now, consider the edge coloring on $\mathbb{Z}_{(n+1)\times (n+1)}$ by using the tiles in $\mathcal{B}$. For $n\geq 5$, let $\mathbb{A}_{n+1}=\mathbb{Z}_{(n+1)\times (n+1)} \setminus \mathbb{Z}_{(n-3)\times (n-3)}((2,2))$ be a rectangularly annular lattice with width 2 (two unit square lattices). The following two steps yield an upper bound of $\Gamma_{(n+1)\times (n+1)}(\mathcal{B})$.

\begin{enumerate}
\item[Step (1):] color the edges on lattice $\mathbb{A}_{n+1}$ by using the four tiles in $\mathcal{B}$. Then, the number of possible patterns is less than $4^{8n-16}$.

\item[Step (2):] color the remaining edges on $\mathbb{Z}_{(n-3)\times (n-3)}((2,2))$ along $\frac{3}{4}\pi$ lines (See Fig. 4.1) by using the tiles in $\mathcal{B}$. Facts (1)$\sim$(3) can be used to prove that the number of possible patterns on each $\frac{3}{4}\pi$ line is less than $4(n-3)$.

\begin{equation*}
\begin{array}{c}
\includegraphics[scale=0.5]{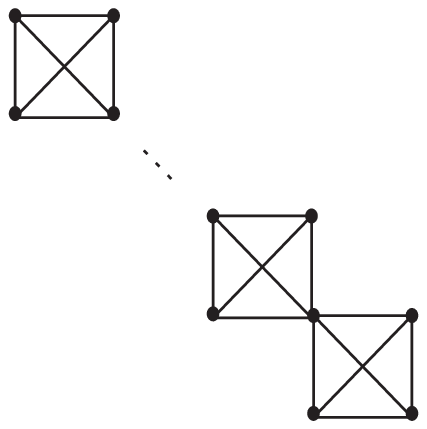}
\end{array}
\end{equation*}
\begin{equation*}
\text{Figure 4.1}
 \end{equation*}
 \end{enumerate}

Accordingly, $\Gamma_{(n+1)\times (n+1)}(\mathcal{B})\leq 4^{8n-16}\cdot(4(n-3))^{2n-9}$.  Therefore, $h(\mathcal{B})=0$. The proof is complete.
$\Box$
\medbreak

Propositions \ref{proposition:4.1}$\sim$\ref{proposition:4.5} suffice to prove the zero entropy of the 31 equivalent classes that are presented in Table A.2.

\begin{theorem}
\label{theorem:4.6}
For the 31 equivalent classes $[\mathcal{B}]$ in Table A.2, if $\mathcal{B}_{M}\in \mathcal{M}(\mathcal{B})$, then $h(\mathcal{B}')=0$ for all $\mathcal{B}'\in [\mathcal{B}_{M}]$. In particular, if $\mathcal{B}''=C_{1}\cup N\subset \Sigma_{2\times 2}(2)$, then $h(\mathcal{B}'')=0$.
\end{theorem}
\textit{Proof}
For each equivalent class $[\mathcal{B}_{M}]$ in Table A.2, there exists $\bar{\mathcal{B}}\in[\mathcal{B}_{M}]$ such that $h(\bar{\mathcal{B}})=0$ can be proven by Proposition \ref{proposition:4.1}$\sim$\ref{proposition:4.5}. Table A.2 presents relevant details. From (\ref{eqn:2.4}), $h(\mathcal{B}')=0$ for all $\mathcal{B}'\in [\mathcal{B}_{M}]$.

The six equivalent classes of $\mathcal{C}(2)$ are denoted with superscript $\ast$ in Table A.2. Therefore, $h(\mathcal{B}'')=0$ for all  $\mathcal{B}''=C_{1}\cup N\subset \Sigma_{2\times 2}(2)$. The proof is complete.
$\Box$
\medbreak

\section{Main results}
\label{5}
In this section, the spatial chaos problem of Wang tiles with two symbols is shown to be determined completely by minimal cycle generators. Furthermore, there are 39 equivalent classes of marginal positive-entropy sets and 18 equivalent classes of saturated zero-entropy sets, which are the "boundary" of between the sets of Wang tiles with positive entropy and the sets of Wang tiles with zero entropy.

First, by Theorems \ref{theorem:3.4} and \ref{theorem:4.6}, for any $\mathcal{B}\subseteq \Sigma_{2\times 2}(2)$, the spatial chaos probelm is determined by the minimal cycle generators in $\mathcal{B}$, as follows.
\begin{theorem}
\label{theorem:4.7}
Given $\mathcal{B}\subset\Sigma_{2\times 2}(2)$, write
$\mathcal{B}=C_{1}\cup C_{2} \cup\cdots \cup C_{k} \cup N$.
Then,
\begin{equation}\label{eqn:4.1}
\begin{array}{ccc}
h(\mathcal{B})>0 & \text{if and only if } & h(C_{1}\cup C_{2} \cup\cdots \cup C_{k})>0.
\end{array}
\end{equation}
\end{theorem}

Next, by carefully studying the 32 sets $\mathcal{B}_{M}$ in Table A.2 and the 39 equivalent classes in Table A.3, the positivity problems of all $\mathcal{B}\subseteq \Sigma_{2\times 2}(2)$ can be solved completely as follows.

\begin{theorem}
\label{theorem:4.7-1}
Consider the sets of Wang tiles with two symbols.
\begin{enumerate}
\item[(i)] There are 39 equivalent classes of marginal positive-entropy sets of Wang tiles that are listed in Table A.3,

\item[(ii)] there are 18 equivalent classes of saturated zero-entropy sets of Wang tiles that are the equivalent classes of $\mathcal{B}_{M}$ in (14)$\sim$(31) in Table A.2.
\end{enumerate}
For a set of Wang tiles $\mathcal{B}$, $h(\mathcal{B})$ is positive if and only if $\mathcal{B}$ contains a marginal positive-entropy set of Wang tiles and $h(\mathcal{B})$ is zero if and only if $\mathcal{B}$ is a subset of a saturated zero-entropy set of Wang tiles.
\end{theorem}
\textit{Proof}
From Theorems \ref{theorem:3.4} and \ref{theorem:4.6}, it is clear that the 39 equivalent classes listed in Table A.3  are the equivalent classes of marginal positive-entropy sets. The detail is omitted here.

Clearly, Theorems \ref{theorem:3.4} and \ref{theorem:4.6} imply that if $\mathcal{B}\subset\Sigma_{2\times 2}(2)$ with $h(\mathcal{B})=0$, then $\mathcal{B}$ is a subset of $\mathcal{B}'\in [\mathcal{B}_{M}]$ for some $\mathcal{B}_{M}$ in Table A.2. Now, the SZE sets can be obtained from the 32 equivalent classes $[\mathcal{B}_{M}]$ in Table A.2 as follows. If $[\mathcal{B}_{1}]$ and $[\mathcal{B}_{2}]$ are two of the 32 equivalent classes $[\mathcal{B}_{M}]$ in Table A.2 with $[\mathcal{B}_{1}]\supseteq [\mathcal{B}_{2}]$, then delete $[\mathcal{B}_{2}]$ from the 32 equivalent classes. It can be proven that the 32 equivalent classes can be reduced to the 18 equivalent classes $[\mathcal{B}_{M}]$ that are the equivalent classes of $\mathcal{B}_{M}$ in (14)$\sim$(31) in Table A.2.

Subsequently, for any set $\mathcal{B}$ in the 18 equivalent classes, it can be verified that if $\mathcal{B}'\supsetneqq \mathcal{B}$, $\mathcal{B}'$ contains a set in the 39 equivalent classes in Table A.3, i.e., $h(\mathcal{B}')>0$. Therefore, the sets in the 18 equivalent classes are the SZE sets.

The proof is complete.
$\Box$
\medbreak

\begin{remark}
\label{remark:4.7-2}
From Table A.3, that the entropies of the sets in the 39 equivalent classes in Table A.3 are positive can be proven by (\ref{eqn:3.8}) with $1\leq m,k\leq 5$. Therefore, for any $\mathcal{B}\subset\Sigma_{2\times 2}(2)$, the positivity of the entropy $h(\mathcal{B})$ can be determined by checking whether or not there exist $1\leq m,k\leq 5$ and $\beta_{j}\in\{1,4\}$, $1\leq j\leq k$, such that $\rho\left( S_{m;\beta_{1}}S_{m;\beta_{2}}\cdots S_{m;\beta_{k}}\right)>1$.
\end{remark}

\begin{remark}
\label{remark:4.8}
Let $\mathcal{B}=C_{1}\cup C_{2} \cup\cdots \cup C_{k} \cup N$. Determining whether $h(\mathcal{B})=h(C_{1}\cup C_{2} \cup\cdots \cup C_{k})$ is difficult, and requires further investigation.

\end{remark}

\begin{remark}
\label{remark:4.9}
In \cite{10}, the decidability of corner coloring with two symbols has been proven, i.e., (1.1) holds. Clearly, in corner coloring, the permutations of symbols in the horizontal and vertical directions are mutually dependent, which means that corner coloring is more closely related in these two perpendicular directions. It can be verified that the positivity of entropy of corner coloring is not determined by minimal cycle generators only.
For example, let
\begin{equation*}
\mathcal{B}_{c;1}=\left\{
\begin{array}{cccc}
\includegraphics[scale=0.7]{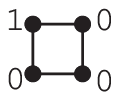},& \includegraphics[scale=0.7]{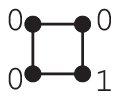},&
\includegraphics[scale=0.7]{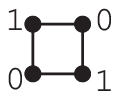}, &\includegraphics[scale=0.7]{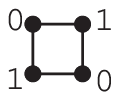}
\end{array}
\right\}
\end{equation*}
and
\begin{equation*}
\mathcal{B}_{c;2}=\mathcal{B}_{c;1}\cup\left\{
\begin{array}{ccccc}
\includegraphics[scale=0.7]{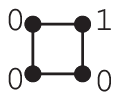},& \includegraphics[scale=0.7]{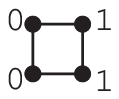},&
\includegraphics[scale=0.7]{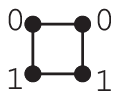}, &\includegraphics[scale=0.7]{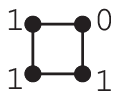}, &\includegraphics[scale=0.7]{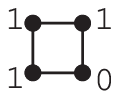}
\end{array}
\right\}.
\end{equation*}
From \cite{10}, $\mathcal{B}_{c;1}$ is a pure set of minimal cycle generators, which are
$\left\{\begin{array}{cc}
\includegraphics[scale=0.7]{ci.eps}, &\includegraphics[scale=0.7]{cj.eps}
\end{array}
 \right\}$ and $\left\{\begin{array}{ccc}
\includegraphics[scale=0.7]{ce1.eps}, &\includegraphics[scale=0.7]{ce4.eps}, &\includegraphics[scale=0.7]{cj.eps}
\end{array}
 \right\}$, and $\mathcal{B}_{c;2}$ contains no minimal cycle generators except these two minimal cycle generators.
As Proposition \ref{proposition:4.1}, it can be easily verified that $\Gamma_{m\times m}(\mathcal{B}_{c;1})\leq 2^{m+n-1}$, and then $h(\mathcal{B}_{c;1})=0$. On the other hand,
by \cite{3}, it can be shown that $h(\mathcal{B}_{c;2})\geq \frac{1}{9}\rho\left( S_{3;11}S_{3;14}S_{3;41}\right)=\frac{1}{9}\log 3 >0$. The details are omitted here.
The complete results of the positivity problem of corner coloring is still under investigation.
\end{remark}

\begin{remark}
\label{remark:4.10}
The concepts of MPE and SZE are also useful in studying the positivity of entropy in other situations. For example, vertex coloring problem and Wang tiles with colors $p\geq 3$.
\end{remark}

\section*{Appendix }

\section*{A.1.} The details of six equivalent classes of $\mathcal{C}(2)$ are
listed in Table A.1.

\begin{equation*}
\begin{tabular}{|lcl|}
\hline
$
[\{O\}]$ & $=$ & $\left\{\{O\},\hspace{0.1cm} \{I\},\hspace{0.1cm}
\{J\},\hspace{0.1cm} \{E\}\right\}
$ \\
\hline
$\left[\{E_{1},E_{4}\}\right]$ & $=$ & $\left\{\{E_{1},E_{4}\},\hspace{0.1cm} \{E_{2},E_{3}\},\hspace{0.1cm}
\{\overline{E}_{1},\overline{E}_{4}\},\hspace{0.1cm} \{\overline{E}_{2},\overline{E}_{3}\}\right\} $ \\
\hline
$\left[\{E_{1},\overline{E}_{1}\}\right]$ & $=$ & $\left\{\{E_{1},\overline{E}_{1}\},\hspace{0.1cm} \{E_{2},\overline{E}_{2}\},\hspace{0.1cm}
\{E_{3},\overline{E}_{3}\},\hspace{0.1cm} \{E_{4},\overline{E}_{4}\}\right\} $ \\
\hline
$\left[\{B,T\}\right]$ & $=$ & $\left\{\{B,T\},\hspace{0.1cm} \{L,R\}\right\} $\\
\hline
$
\left[\{E_{1},B,R\}\right]$ & $=$ & $\left\{
\begin{array}{cccc}\{E_{1},B,R\}, & \{E_{2},B,L\},&
\{E_{3},T,R\},&\{E_{4},T,L\},\\
 \{\overline{E}_{1},T,L\},& \{\overline{E}_{2},T,R\},&
\{\overline{E}_{3},B,L\},& \{\overline{E}_{4},B,R\}
\end{array}\right\}
$ \\
\hline
$
\left[\{E_{1},E_{2},B\}\right]$& $=$ & $\left\{
\begin{array}{cccc}\{E_{1},E_{2},B\},&\{E_{1},E_{3},R\},&
\{E_{2},E_{4},L\},& \{E_{3},E_{4},T\},\\
 \{E_{1},\overline{E}_{2},R\},& \{E_{1},\overline{E}_{3},B\},&
\{E_{2},\overline{E}_{1},L\},& \{E_{2},\overline{E}_{4},B\} , \\
 \{E_{3},\overline{E}_{1},T\},& \{E_{3},\overline{E}_{4},R\},&
\{E_{4},\overline{E}_{2},T\},& \{E_{4},\overline{E}_{3},L\},\\
\{\overline{E}_{1},\overline{E}_{2},T\},&\{\overline{E}_{1},\overline{E}_{3},L\},&
\{\overline{E}_{2},\overline{E}_{4},R\},& \{\overline{E}_{3},\overline{E}_{4},B\}
\end{array}\right\}
$ \\
\hline
\end{tabular}%
\end{equation*}

\begin{equation*}
\text{Table A.1}
\end{equation*}

\section*{A.2.}
For $1\leq k\leq 5$, Table A.2 presents the 31 equivalent classes of $\mathcal{B}=C_{1}\cup C_{2} \cup\cdots \cup C_{k}\subset \Sigma_{2\times 2}(2)$, for which $h(\mathcal{B})>0$ can not be proven by Theorem \ref{theorem:3.1}. Table A.2 also presents all $\mathcal{B}_{M}$ of the 31 equivalent classes and the corresponding propositions that are used to prove $h([\mathcal{B}_{M}])=0$.

For simplicity, every tile $B\in\Sigma_{2\times 2}(2)$ with $(b(B),l(B),t(B),r(B))=(u_{1},u_{2},u_{3},u_{4})$ is assigned a number according to $\psi(B)\equiv\psi(u_{1},u_{2},u_{3},u_{4})$.
Indeed,
\begin{equation*}
\left[
\begin{array}{cccc}
O & E_{2} & E_{4} & R \\
E_{3} & J & B & \overline{E}_{1} \\
E_{1} & T & I & \overline{E}_{3}\\
L &  \overline{E}_{4} &  \overline{E}_{2} & E
\end{array}
\right]=
\left[
\begin{array}{cccc}
1 & 3 & 2 & 4 \\
9 & 11 & 10 & 12 \\
5 & 7 & 6 & 8 \\
13 & 15 & 14 & 16
\end{array}
\right].
\end{equation*}

{\small
\begin{equation*}
\begin{tabular}{|l|c|c|}
  \hline
   $[\mathcal{B}]=[C_{1}\cup C_{2} \cup\cdots \cup C_{k}]$  &  $N=\mathcal{B}_{M}\setminus \mathcal{B}$  & Proposition used for $h\left([\mathcal{B}_{M}]\right)=0$  \\ \hline
1*. \hspace{0.1cm}$\left[\{1\}\right]$   &  \{2, 3, 4, 7, 8, 12\} &  4.2   \\ \hline
2. \hspace{0.1cm}$\left[\{1,6 \}\right]$   &  \{2, 3, 4, 7, 8, 12\} & 4.2    \\ \hline
   & \{5, 7, 9, 13, 14, 15\} &  4.4 \\ \hline
3.  \hspace{0.1cm}$\left[ \{1, 16 \}\right]$ &  \{2, 3, 4, 7, 8, 12\} &4.2 \\ \hline
4*. \hspace{0.1cm}$\left[\{2,5\}\right]$   &  \{9, 10, 13, 14\} & 4.4    \\ \hline
5*. \hspace{0.1cm}$\left[\{5,12\}\right]$   &  \{9, 10, 13, 14\} &  4.4   \\ \hline
6*. \hspace{0.1cm}$\left[\{7,10\}\right]$   &  \{2, 3, 4, 8, 12\} & 4.3    \\ \hline
7.  \hspace{0.1cm}$\left[\{1, 5, 12 \}\right]$&  \{9, 10, 13, 14\} & 4.4 \\ \hline
8.  \hspace{0.1cm}$\left[\{1, 6, 11\} \right]$&  \{2, 3, 4, 7, 8, 12\} & 4.2\\\hline
9.  \hspace{0.1cm}$\left[\{1, 7, 10 \}\right]$&  \{2, 3, 4, 8, 12\} &  4.3\\\hline
10.  \hspace{0.1cm}$\left[\{1, 12, 15\}\right] $& \{9, 10, 13, 14\} & 4.4 \\\hline
11.  \hspace{0.1cm}$\left[\{2, 5, 12\} \right]$&  \{9, 10, 13, 14\} & 4.4\\\hline
12*. \hspace{0.1cm}$\left[\{3,5,10\}\right]$   &  $\emptyset$ & 4.1    \\ \hline
13*. \hspace{0.1cm}$\left[\{4,5,10\}\right]$   &  $\emptyset$ &  4.1   \\ \hline
14.  \hspace{0.1cm}$\left[\{1, 4, 5, 14\}\right]$ &   $\emptyset$ &  4.5 \\\hline
15.  \hspace{0.1cm}$ \left[\{1, 5, 12, 15 \}\right]$&  \{9, 10, 13, 14\} & 4.4 \\\hline
16.  \hspace{0.1cm}$\left[ \{1, 5, 12, 16\} \right]$&  \{9, 10, 13, 14\} & 4.4\\\hline
17.  \hspace{0.1cm}$\left[\{1, 6, 11, 16\} \right]$&  \{2, 3, 4, 7, 8, 12\} &  4.2\\\hline
18.  \hspace{0.1cm}$\left[\{1, 6, 12, 15 \} \right]$& \{9 ,10, 13, 14\} &  4.4\\\hline
19.  \hspace{0.1cm}$\left[\{1, 7, 10, 16 \} \right]$&\{2, 3, 4, 8, 12\} &  4.3\\\hline
20.  \hspace{0.1cm}$\left[\{1, 7, 12, 13\}\right]$&   $\emptyset$ &  4.1\\\hline
21.  \hspace{0.1cm}$\left[ \{1, 7, 12, 14\} \right]$&   $\emptyset$ & 4.1 \\\hline
22.  \hspace{0.1cm}$ \left[\{2, 5, 12, 15 \}\right]$&  \{9, 10, 13, 14\} & 4.4 \\\hline
23.  \hspace{0.1cm}$\left[ \{3, 4, 5, 10 \}\right]$&   $\emptyset$ &  4.1\\\hline
24.  \hspace{0.1cm}$\left[\{3, 5, 7, 10\} \right]$&   $\emptyset$ & 4.1\\\hline
25.  \hspace{0.1cm}$\left[\{3, 5, 8, 10 \}\right]$& $\emptyset$ & 4.1\\\hline
26.  \hspace{0.1cm}$\left[ \{3, 5, 10, 12\}\right]$ & $\emptyset$ & 4.1 \\\hline
27.  \hspace{0.1cm}$ \left[\{3, 5, 12, 14 \}\right]$&  $\emptyset$ & 4.1\\\hline
28.  \hspace{0.1cm}$\left[ \{4, 5, 7, 10\}\right]$&    $\emptyset$ & 4.1 \\\hline
29.  \hspace{0.1cm}$\left[\{4, 5, 10, 15 \}\right]$&   $\emptyset$ & 4.1 \\\hline
30.  \hspace{0.1cm}$\left[ \{4, 7, 10, 13\}\right]$&   $\emptyset$ &4.1\\\hline
31.  \hspace{0.1cm}$ \left[\{5, 7, 10, 12\}\right]$ &  $\emptyset$ &4.1\\ \hline
\end{tabular}
\end{equation*}
}

\begin{equation*}
\text{Table A.2}
\end{equation*}

\section*{A.3.}
Table A.3 presents the 39 equivalent classes that are introduced in the proof of Theorem \ref{theorem:3.4}. By (\ref{eqn:3.8}), it can be verified that the 39 equivalent classes have positive entropies. For each $[\mathcal{B}]$ of the 39 equivalent classes, Table A.3 also presents the values $m$ and $\beta_{1},\beta_{2}, \cdots,\beta_{k}$ of (\ref{eqn:3.8}) for  proving $h(\mathcal{B})>0$.

{\small
\begin{equation*}
\begin{tabular}{|l|c|c|}
  \hline
   $[\mathcal{B}]=[C_{1}\cup C_{2} \cup\cdots \cup C_{k}]$  &  $m$ & $\beta_{1}, \beta_{2}, \cdots, \beta_{k}$ \\ \hline
1. \hspace{0.1cm}$\left[\{1, 2, 5\}\right]$   &  2 &  1   \\ \hline
2. \hspace{0.1cm}$\left[\{1, 3, 5, 10 \}\right]$    &  2 &  1   \\ \hline
3.  \hspace{0.1cm}$\left[ \{1, 3, 6, 14 \}\right]$  &  1 &  1, 4   \\ \hline
4. \hspace{0.1cm}$\left[\{1, 4, 5, 10\}\right]$   &  2 &  1  \\ \hline
5. \hspace{0.1cm}$\left[\{1, 5, 6, 12\}\right]$  &  3 &  1   \\ \hline
6. \hspace{0.1cm}$\left[\{1, 6, 7, 10\}\right]$    &  3 &  1, 4   \\ \hline
7.  \hspace{0.1cm}$\left[\{2, 3, 5, 9 \}\right]$ &  2 &  1  \\ \hline
8.  \hspace{0.1cm}$\left[\{2, 3, 5, 10\} \right]$ &  5 &  1, 1, 4   \\ \hline
9.  \hspace{0.1cm}$\left[\{2, 3, 5, 14 \}\right]$ &  2 &  1, 4   \\ \hline
10.  \hspace{0.1cm}$\left[\{2, 4, 5, 10\}\right] $ &  2 &  1   \\ \hline
11.  \hspace{0.1cm}$\left[\{2, 5, 7, 10\} \right]$ &  2 &  1, 4   \\ \hline
12. \hspace{0.1cm}$\left[\{4, 5, 10, 12\}\right]$    &  2 &  1   \\ \hline
13. \hspace{0.1cm}$\left[\{1, 3, 5, 12, 13\}\right]$    &  3 &  1, 1, 1, 4   \\ \hline
14.  \hspace{0.1cm}$\left[\{1, 3, 5, 12, 14\}\right]$  &  4 &  1, 1, 1, 1, 4   \\ \hline
15.  \hspace{0.1cm}$ \left[\{1, 3, 6, 12, 13 \}\right]$ &  4 &  1, 1, 4   \\ \hline
16.  \hspace{0.1cm}$\left[ \{1, 3, 7, 12, 13\} \right]$ &  3 &  1, 1, 1, 1, 4   \\ \hline
17.  \hspace{0.1cm}$\left[\{1, 4, 5, 12, 14\} \right]$ &  3 &  1, 1, 4, 4   \\ \hline
18.  \hspace{0.1cm}$\left[\{1, 4, 5, 14, 15 \} \right]$ &  4 &  1, 4   \\ \hline
19.  \hspace{0.1cm}$\left[\{1, 4, 5, 14, 16 \} \right]$ &  4 &  1, 1, 4, 4   \\ \hline
20.  \hspace{0.1cm}$\left[\{1, 4, 7, 10, 13\}\right]$ &  2 &  1   \\ \hline
21.  \hspace{0.1cm}$\left[ \{1, 5, 7, 8, 10\} \right]$ &  3 &  1, 1, 4   \\ \hline
22.  \hspace{0.1cm}$ \left[\{1, 5, 7, 10, 12 \}\right]$ &  3 &  1, 1, 1, 4   \\ \hline
23.  \hspace{0.1cm}$\left[ \{1, 5, 7, 12, 14 \}\right]$ &  3 &  1   \\ \hline
24.  \hspace{0.1cm}$\left[\{1, 6, 7, 12, 13\} \right]$ &  3 &  1   \\ \hline
25.  \hspace{0.1cm}$\left[\{1, 7, 10, 12, 13 \}\right]$ &  3 &  1, 1, 1, 1, 4   \\ \hline
26.  \hspace{0.1cm}$\left[ \{1, 7, 10, 12, 14\}\right]$  &  3 &  1, 1, 4   \\ \hline
27.  \hspace{0.1cm}$ \left[\{1, 7, 12, 13, 14 \}\right]$ &  4 &  1  \\ \hline
28.  \hspace{0.1cm}$\left[ \{3, 4, 5, 7, 10\}\right]$ &  5 &  1, 1, 4   \\ \hline
29.  \hspace{0.1cm}$\left[\{3, 4, 5, 8, 10\}\right]$ &  5 &  1   \\ \hline
30.  \hspace{0.1cm}$\left[ \{3, 4, 5, 10, 15\}\right]$ &  2 &  1, 1, 4   \\ \hline
\end{tabular}
\end{equation*}
}

{\small
\begin{equation*}
\begin{tabular}{|l|c|c|}
  \hline
31.  \hspace{0.1cm}$ \left[\{3, 5, 7, 8, 10\}\right]$ &  4 &  1, 4   \\ \hline
32. \hspace{0.1cm}$\left[\{3, 5, 7, 10, 12\}\right]$    &  5 &  1, 1, 4   \\ \hline
33. \hspace{0.1cm}$\left[\{3, 5, 8, 10, 12\}\right]$   &  4 &  1, 1, 1, 4, 4   \\ \hline
34.  \hspace{0.1cm}$\left[ \{3, 5, 8, 10, 15 \}\right]$ &  4&  1, 4   \\ \hline
35. \hspace{0.1cm}$\left[\{3, 5, 10, 12, 13\}\right]$   &  2 &  1  \\ \hline
36. \hspace{0.1cm}$\left[\{3, 5, 10, 12, 14\}\right]$    &  3 &  1, 1, 4, 4   \\ \hline
37. \hspace{0.1cm}$\left[\{4, 5, 7, 9, 10\}\right]$   &  2 &  1   \\ \hline
38.  \hspace{0.1cm}$\left[\{4, 5, 7, 10, 13 \}\right]$ &  2 &  1   \\ \hline
39.  \hspace{0.1cm}$\left[\{4, 5, 7, 10, 15\} \right]$ &  2 &  1, 4   \\ \hline
\end{tabular}
\end{equation*}

\begin{equation*}
\text{Table A.3}
\end{equation*}







\end{document}